\documentclass{amsart}
\usepackage{amsfonts}
\usepackage{amsmath,amsthm}
\usepackage{amsfonts,amssymb}

\usepackage{enumerate}

\newtheorem{thm}{Theorem}[section]

\newtheorem{lem}[thm]{Lemma}

\newtheorem{defn}[thm]{Definition}
\newtheorem{rem}[thm]{Remark}

\numberwithin{equation}{section}

\newcommand{\al}{\alpha}

\def\sph{\mathbb{S}^{d}}

\def\l{\left}

\def\vz{\varepsilon}

\def\lz{\lambda}

\def\az{\alpha}

\def\sz{\sigma}

\def\f{\frac}
\def\({\Bigl(}
\def \){ \Bigr)}

\def\Ga{\Gamma}

\def\sub{\substack}

 \def\l{{\lambda}}

 \def\rb{{\mathbf r}}

 \def\RR{{\mathbb R}}

\def\sa{\sigma}
\def\sz{\sigma}

\def\ss{{\Bbb S}^{d}}

\def\lb{\langle}
\def\rb{\rangle}

\def\bi{\bibitem}

\begin{document}
\def\RR{\mathbb{R}}
\def\Exp{\text{Exp}}
\def\FF{\mathcal{F}_\al}

\title[] {Approximation and quadrature  by weighted least squares polynomials on the sphere}

\author{Wanting Lu} \address{ School of Mathematical Sciences, Capital Normal
University, Beijing 100048,
 China.}
 \email{luwanting1234@163.com}

\author{Heping Wang} \address{ School of Mathematical Sciences, Capital Normal
University,
Beijing 100048,
 China.}
\email{wanghp@cnu.edu.cn.}

\keywords{Marcinkiewicz-Zygmund inequalities; Sampling; Sobolev
spaces; Least squares problem;  operator norms}
\subjclass[2010]{41A17; 42A10; 65D30.}

\begin{abstract} Given a sequence of Marcinkiewicz-Zygmund inequalities in $L_2$ on  a usual compact space $\mathcal M$,
Gr\"ochenig in \cite{Gr} introduced the weighted least squares
polynomials and the least squares quadrature from  pointwise
samples of a function, and  obtained approximation theorems and
quadrature errors.  In this paper
  we
  obtain approximation theorems and quadrature errors on the sphere which are optimal. We also give  upper bounds of  the operator norms of the
weighted least squares operators.

\end{abstract}

\maketitle
\input amssym.def

\section{Introduction}

Let $$\ss:=\{x=(x_1,\dots,x_{d+1})\in\mathbb{R}^{d+1}\,\big|\
|x|^2=\sum_{k=1}^{d+1}|x_k|^2=1\}\ \ \ (d\ge 2)
$$ be the unit sphere in $\mathbb{R}^{d+1}$ endowed with the
rotationally  invariant measure $\sz$ normalized by
$\int_{\ss}d\sz=1$. We denote by $L_2(\ss)$ the usual Hilbert
space of square-integrable functions on $\ss$ with the inner
product\begin{equation} \label{1.0}\langle
f,g\rangle:=\int_{\ss}f(x)g(x)d\sz(x)\end{equation} and the norm
$\|f\|_2=\langle f,f\rangle^{1/2}$, and by $C(\ss)$ the space of
continuous functions on $\ss$ with supremum norm
$$\|f\|:=\|f\|_\infty:=\sup_{x\in\ss}|f(x)|.$$

The space $\Pi_n^d$   of  spherical polynomials on $\ss$ of degree
 $\le n$ consists of the restrictions to $\ss$ of all
polynomials on $\Bbb R^{d+1}$ of total degree $\le n$. The
dimension of $\Pi_n^d$ is  given dy
$$d_n:= \dim(\Pi_n^d)=\frac{(2n+d)\Gamma(n+d)}{\Gamma(d+1)\Gamma(n+1)}.$$

Let $S_{n}$ be the orthogonal projection from $ L_{2}(\ss)$ onto
$\Pi_{n}^{d}$, i.\! e., for $f\in L_2(\ss)$,
\begin{equation}\label{1.1}S_{n}f({x})=\int_{\ss}f( y)E_{n}(x,y)d\sz(y),\end{equation}where $E_{n}(x,y)$ is the
reproducing kernel for $\Pi_{n}^{d}$ with respect to the inner
product \eqref{1.0}. We note that $E_{n}(x,y)$ satisfies the
following properties:

(1) For any $x,y\in \ss$, $E_{n}(x,y)=E_{n}(y,x)$;

(2) For any fixed $y\in \ss$, $E_{n}(\cdot,y)\in\Pi_{n}^{d}$;

(3) For any $ p\in\Pi_{n}^{d}$ and $x\in \ss$, $p(x)=S_np(x)=
\langle p,E_{n}(\cdot,x)\rangle$.

This paper is concerned with constructive polynomial approximation
on $\ss$ that uses function values (the samples) at selected
well-distributed points (sometimes called standard information).
Here the ``well-distributed'' points indicate that those points
constitute a  Marcinkiewicz-Zygmund (MZ) family on $\ss$ defined
as follows.

\begin{defn}Let $\mathcal X  =\{\mathcal X_n\}= \{x_{n,k} : n=1,2,\dots, \, k = 1,\dots,
l_n\}$ be a doubly-indexed set of points in $\ss$ and $\tau=
\{\tau_{n,k} : n=1,2,\dots, \, k = 1,\dots, l_n\}\subset
(0,+\infty)$ be a family of non-negative weights. Then $\mathcal X
$ is called a Marcinkiewicz-Zygmund (MZ) family with associated
weight $\tau$, if there exist constants $A, B
> 0$ independent of $n$ such that \begin{equation}\label{1.2}A\|p\|_2^2\le \sum_{k=1}^{l_n}|p(x_{n,k})|^2\tau_{n,k}\le B\|p\|_2^2\ \qquad {\rm
for\ all}\ p\in\Pi_n^d. \end{equation} The ratio $\kappa = B/A$ is
the global condition number of the Marcinkiewicz-Zygmund family
$\mathcal X$, and $\mathcal X _n = \{x_{n,k} : k = 1, \dots,
l_n\}$ is the $n$th layer of $\mathcal X$.\end{defn}

The MZ inequality \eqref{1.2} means that the $L_{2}$-norm  of a
spherical polynomial of degree at most $n$ is comparable to the
discrete version given by the weighted $\ell_2$-norm of its
restriction to $\mathcal X_n$. It follows from \cite{MNW, BD,
DX13} that such MZ families exist if the families are dense
enough. We remark that Fekete points of degree $\lfloor
n(1+\vz)\rfloor\ (\vz>0)$ on the sphere are  MZ families with the
equal weights $\tau_{n,k}=\frac 1{d_n}$, see \cite{MO}, where
$\lfloor a\rfloor$ is the largest integer not exceeding $a\in\Bbb
R$. Also, sufficient conditions and necessary density conditions
for MZ families with the  equal weights $\tau_{n,k}=\frac 1{d_n}$
on the sphere are obtained in \cite{MP} and \cite{M},
respectively.

Given   a MZ family on a usual compact space $\mathcal M$,
Gr\"ochenig in \cite{Gr} introduced the weighted least squares
polynomials and the least squares quadrature from pointwise
samples of a function, and  obtained approximation theorems and
quadrature errors. However, the obtained error estimates in
\cite{Gr} are not optimal due to the generality of $\mathcal M$.
In this paper we confine the study to the sphere $\ss$.

 Let
$\mathcal X$ be a MZ family with associated weight $\tau$. Given
the samples $\{ f (x_{n,k} )\}$ of a continuous function $f$ on
$\ss$ on the $n$th layer $\mathcal X_n$, we want to approximate
$f$ using only these samples. For this we solve a sequence of
weighted least squares problems with samples taken from the
samples $\mathcal X_n$:
\begin{equation}\label{1.3}L_n f=\arg \min_{p\in\Pi_n^d}\
\sum_{k=1}^{l_n}|f(x_{n,k})-p(x_{n,k})|^2\tau_{n,k}.\end{equation}
This procedure yields  a sequence $\{L_nf\}$ of the best weighted
$\ell_2$-approximation of the data $\{ f (x_{n,k} )\}$ by   a
  spherical polynomial in $\Pi_n^d$ for every $f\in
C(\ss)$. We call  $L_nf$   the weighted least squares polynomial,
and $L_n$   the weighted least squares operator.
 Clearly, $L_n $
is the projection onto $\Pi_{n}^{d}$, i.e., $L_n $ is a bounded
linear operator on $C(\ss)$ satisfying that $L_n^2=L_n $ and the
range of $L_n $ is $\Pi_{n}^{d}$. It follows from \eqref{1.2} that
for $p\in \Pi_n^d$, $p=0$ whenever $p(x_{n,k})=0$. This means that
usually $\mathcal X_n$ contains more than $d_n=\dim \Pi_n^d$
points, so that it is not an interpolating set for $\Pi_n^d$.
Therefore, $L_nf$ is usually  a quasi-interpolant.

 Let
$E_n(x,y)$ be the reproducing kernel of $\Pi_n^d$ with respect to
the inner product \eqref{1.0}. Then the Marcinkiewicz-Zygmund
inequalities \eqref{1.2} say that every set $\{\tilde e_{n,k}: k =
1,\dots , l_n\},\ \tilde e_{n,k}=\tau^{ 1/2}_{
n,k}E_n(\cdot,x_{n,k}) $ is a frame for $\Pi_n^d$ with uniform
frame bounds $A, B
> 0$, i.e., for all $p\in\Pi_n^d$, we have
$$A \|p\|_2^2 \le \sum_{k=1}^{l_n}|\langle p, \tilde e_{n,k}\rangle|^2\le
B\|p\|_2^2.$$ It follows that the associated frame operator $$T_n
p(\cdot) = \sum_{k=1}^{l_n} \langle p, \tilde e_{n,k}\rangle
\tilde
e_{n,k}(\cdot)=\sum_{k=1}^{l_n}\tau_{n,k}p(x_{n,k})E_n(\cdot,x_{n,k})$$
is invertible on $\Pi_n^d$  for every $n\in \Bbb N$. We obtain the
dual frame $ \{e_{n,k} = T_n^{-1} (\tilde e_{n,k}), \
k=1,\dots,l_n\}$ for $\Pi_n^d$  with uniform frame bounds $B^{-1}$
and $A^{-1}$, and every polynomial $p\in\Pi_n^d$ can be
reconstructed from the samples on $\mathcal X_n$  by
$$p=T_n^{-1}T_np=\sum_{k=1}^{l_n} \langle p, \tilde e_{n,k}\rangle T_n^{-1}\tilde e_{n,k}=\sum_{k=1}^{l_n} \tau_{n,k}^{1/2}p(x_{n,k})e_{n,k}. $$
The weights for the quadrature rules are defined by
$$w_{n,k} = \langle\tau^{ 1/2}_{ n,k}e_{n,k} , 1\rangle = \tau^{
1/2}_{ n,k}\int_{\ss} e_{n,k} (x) d\sz(x)$$ and the corresponding
quadrature rule is defined by
$$I_n(f)=\sum_{k=1}^{l_n}w_{n,k}f(x_{n,k}).$$
Such quadrature rules are usually named least squares quadrature.

 Let
$\mathcal X$ be a MZ family with associated weight $\tau$. For  a
function $f$ in the Sobolev space $H^\sz(\mathbb{S}^d),\ \sz>d/2$
(see Section 2 for definition of  $H^\sz(\ss)$), Gr\"ochenig
obtained in \cite{Gr} that
$$\|f-L_n f\|_2\le c(1+\kappa^2)^{1/2}n^{-\sz+d/2}\|f\|_{H^\sz},$$
and $$\big|\int_{\ss}f(x)d\sz(x)-I_n (f)\big|\le
c(1+\kappa^2)^{1/2}n^{-\sz+d/2}\|f\|_{H^\sz},$$ where $c$ depends
on $ d,\ s$,  but not on $f$ or $\kappa$ or the MZ family
$\mathcal X$.

We remark that the condition $\sz>d/2$ is the well known necessary
and sufficient condition for $H^\sz(\mathbb{S}^d)$ to be
continuously embedded in $C(\mathbb{S}^d)$, and is unavoidable in
any approximation scheme that employs function values.

In this paper we improve the above results and obtain the optimal
estimates. One of our main results can be formulated as follows.

\begin{thm}Let $\mathcal X$ be a MZ family on $\ss$  with associated weight $\tau$ and global condition
number $\kappa = B/A$, $L_n$ be the weighted least squares
operators defined by \eqref{1.3}, and let $\{I_n \,|\, n\in\Bbb
N\}$ be the sequence of the least squares quadrature. If $f\in
H^{\sz}(\ss),\ \sz>d/2$, then we have
\begin{equation}\label{1.4}\|f-L_n f\|_2\le
c(1+\kappa^2)^{1/2}n^{-\sz}\|f\|_{H^\sz},\end{equation} and
\begin{equation}\label{1.5}\Big|\int_{\ss}f(x)d\sz(x)-I_n (f)\Big|\le
c(1+\kappa^2)^{1/2}n^{-\sz}\|f\|_{H^\sz},\end{equation} where $c$
depends on $ d,\ s$,  but not on $f$ or $\kappa$ or the MZ family
$\mathcal X$.
\end{thm}

Remarkably,  the convergence rates for the weighted least squares
approximation and for the least squares quadrature errors are
optimal (as explained in Remark 1.3 below)  in a variety of
Sobolev space settings. When the global condition number $\kappa =
1$,  the weighted least squares operators are  just the
hyperinterpolation operators  on the sphere (see Subsection 2.5).
In this case,  the inequality \eqref{1.4} was obtained in
\cite{Wang}.

\begin{rem}For $N\in \Bbb N$, the sampling numbers (or the  optimal recovery) of the Sobolev classes
$BH^\sz(\ss)$ (the unit ball of the space $H^\sz(\ss)$) in
$L_2(\ss)$ are defined by
$$g_N(BH^\sz(\ss),L_2(\ss)):=\inf_{\sub{\xi_1,\dots,\xi_N\in\ss\\ \varphi:\,\Bbb R^N\to L_2(\ss)}}
\sup_{f\in BH^\sz(\ss)}\|f-\varphi(f(\xi_1),\dots,f(\xi_N))\|_2,
$$ where the infimum is taken over all $N$ points
$\xi_1,\dots,\xi_N$ in $\ss$ and all mappings $\varphi$ from $\Bbb
R^N$ to $L_2(\ss)$. The optimal quadrature errors of the Sobolev
classes $BH^\sz(\ss)$ are defined by
\begin{equation*}
e_N ( \text{Int}; BH^\sz(\ss)):= \inf_{\sub{\l_1, \dots, \l_N \in
\RR\\
\xi_1,\dots, \xi_N \in \sph}} \sup_{f\in BH^\sz(\ss)}
\Bigl|\int_{\sph} f(x) \, d\sa(x)-\sum_{j=1}^N \l_j
f(\xi_j)\Bigr|.\end{equation*}

It follows from \cite{WW}, \cite{BH} and \cite{H} that
\begin{equation}\label{1.6}
g_N(BH^\sz(\ss),L_2(\ss))\asymp N^{-\sz/d}\ \ \ {\rm and}\ \ \ e_N
( \text{Int}; BH^\sz(\ss))\asymp N^{-\sz/d},
\end{equation}
where the notation $A(N)\asymp B(N) $ means that $A(N)\ll B(N)$
and $A(N)\gg B(N)$,  $A(N)\ll B(N)$ means that there exists a
positive constant $c$ independent of $N$ such that $A(N)\leq
cB(N)$, and $A(N)\gg B(N)$ means $B(N)\ll A(N)$.

 It follows from \cite{M, MO, MP} that there exist MZ families with $l_n\asymp d_n\asymp N\asymp n^d$.   For such MZ family it follows from
\eqref{1.6} that
$$\sup_{f\in BH^\sz(\ss)}\|f-L_n(f)\|_2\asymp N^{-s/d}\asymp
g_N(BH^\sz(\ss),L_2(\ss)),$$which implies that  the weighted least
squares operators $L_n$ are asymptotically optimal algorithms in
the sense of optimal recovery. Also, we have $$\sup_{f\in
BH^\sz(\ss)} \Bigl|\int_{\sph} f(x) \, d\sa(x)-I_n(f)\Bigr|\asymp
N^{-s/d}\asymp  e_N ( \text{Int}; BH^\sz(\ss)),
$$ which means that  the least squares quadrature rules are the
asymptotically optimal quadrature  formulas.
\end{rem}

 For a linear operator $L$ on $C(\ss)$, the
operator norm $\|L\|$ of $L$ is given by
$$\|L\|:=\sup\{\|Lf\|\ \big|\ f\in C(\ss),\ \|f\|\le 1\}.$$
The quantity $\|L\|$ of a projection $L$ is also called the
Lebesgue constant of $L$ and the estimation of $\|L\|$ is
extremely important in numerical computation. We use the
Christoffel functions to get the following estimation.

\begin{thm} Let $\mathcal X$ be a MZ family on $\ss$  with associated weight $\tau$ and global condition
number $\kappa = B/A$,   and let $L_n$ be the weighted least
squares operators defined by \eqref{1.3}. Then we have
$$n^{(d-1)/2}\ll \|L_n\|\ll\kappa^{1/2} n^{d/2}.
$$
\end{thm}

\begin{rem} When $\kappa=1$, the weighted least
squares operators are  just the hyperinterpolation operators on
the sphere (see Subsection 2.5). In this case, we have (see
\cite{SW0, GS, R})
$$\|L_n\|\asymp n^{(d-1)/2}.$$  Hence, we conjecture that the
upper bound of the operator norm of $L_n$  is
 $n^{(d-1)/2}$ multiplied by a positive constant independent of $n$.  However, we have not been able to prove it.\end{rem}

 \begin{rem} Since the weighted least
squares operators  $L_{n}$ is the projection onto $\Pi_{n}^{d}$,
by the Lebesgue theorem the hyperinterpolation error in the
uniform norm can be estimated as
$$\|f-L_n(f)\|\le (1+\|L_n\|)E_n(f)\ll  \kappa^{1/2} n^{d/2}E_n(f),$$ where
$$E_n(f):=\inf_{p\in\Pi_n^d}\|f-p\|$$ is the best approximation of
$f$ by $\Pi_n^d$.
\end{rem}

\section{Preliminary}

This section is devoted to giving some basic facts about spherical
harmonics (see for example, \cite{DX13}) and some preliminary
knowledge.

\subsection{Harmonic analysis on the sphere}\

Let $\ss=\{x\in\Bbb R^{d+1}\ |\ |x|=1\}$ denote the unit sphere in
$\Bbb R^{d+1}$,   where $(x,y)=x\cdot y$ is the usual inner
product and $|x|=(x, x)^{1/2}$ is the Euclidean norm.
 We denote by $\mathcal{H}_\ell^d$ the space of
all spherical harmonics of degree $l$ on $\ss$, i.e., the space of
the restrictions to $\ss$ of all harmonic homogeneous polynomials
of exact degree $\ell$ on $\Bbb R^{d+1}$. It is well known that
the dimension of $\mathcal{H}_\ell^d$ is
$$ N(d,\ell)={\rm dim}\,\mathcal{H}_\ell^d=\left\{\begin{array}{cl} 1,\ \
\
   & {\rm if}\ \ \ell=0,\\
\frac{(2\ell +d-1)\,(\ell +d-2)!}{(d-1)!\ l!},\ \     & {\rm if}\
\ \ell=1,2,\dots.
\end{array}\right.$$  The
spaces $\mathcal{H}_k^d, \ k=0,1,2, \dots$ are just the
eigenspaces corresponding to the eigenvalues $-k(k+d-1)$ of the
Laplace-Beltrami operator $\triangle$ on the sphere $\ss$ and  are
mutually orthogonal with respect to the inner product $$\langle
f,g\rangle =\int_{\ss}f(x)g(x)d\sz(x).
$$

Let $$\{Y_{\ell,k}\equiv Y_{\ell,k}^{d}\ |\ k=1,\dots,
N(d,\ell)\}$$ be a fixed orthonormal basis for
$\mathcal{H}_\ell^d$. We have the addition theorem for the
spherical harmonics of degree $\ell$:
$$\sum_{k=1}^{N(d,\ell)}Y_{\ell,k}(x)Y_{\ell,k}(y)=\frac{\ell
+\lambda}{\lambda}C_\ell^{\lambda}( x\cdot y),\  \ x,y\in\ss,\
\lambda=\frac{d-1}{2},\ \ d\ge2,$$where $C_\ell^\lambda (t)$ is
the usual ultraspherical (Gegenbeuer) polynomial of order
${\lambda}$ normalized by $C_\ell^{\lambda}(1)=\frac{\Gamma(\ell
+2\lambda)}{\Gamma(2\lambda)\Gamma(\ell +1)}$ and  generated by
$$ \frac 1{(1-2zt+z^2)^{\lambda}}=\sum_{k=0}^\infty
C_k^\lambda (t)z^k \ \  (0\leq z <1).$$ (See (\cite[p. 81]{Sz}).

 Since the space $\Pi_n^d$ can be expressed as a direct sum of
\begin{equation*}\Pi_n^d=\mathcal{H}_0^d
\bigoplus\mathcal{H}_1^d\bigoplus\cdots\bigoplus\mathcal{H}_n^d
,\end{equation*}we obtain that $$\{Y_{\ell,k}\ |\ k=1,\dots,
N(d,\ell),\ \ell=0,1,\dots, n\}$$ forms an $L_2$-orthonormal basis
of $\Pi_n^d$. Hence, the $L_2$-reproducing kernel $E_n(x,y)$ of
$\Pi_{n}^{d}$ has the explicit expression:
\begin{equation*}\label{2.00}E_n(x,y)=\sum_{\ell=0}^n\sum_{k=1}^{N(d,\ell)}Y_{\ell,k}(x)Y_{\ell,k}(y)=\sum_{\ell=0}^n\frac{\ell
+\lambda}{\lambda}C_\ell^{\lambda}( x\cdot y),  \
\lz=\frac{d-1}2,\  x,y\in\ss.\end{equation*} Specifically, we have
\begin{equation}\label{2.01}E_n(x,x)=\sum_{\ell=0}^n\frac{\ell
+\lambda}{\lambda}\frac{\Gamma(\ell
+2\lambda)}{\Gamma(2\lambda)\Gamma(\ell +1)}\asymp n^{d},\ \
\lz=\frac{d-1}2, \ \ x\in\ss.\end{equation}

We remark that $$\{Y_{\ell,k}\ |\ k=1,\dots, N(d,\ell),\
\ell=0,1,2,\dots\}$$ is an orthonormal basis for the Hilbert space
$L_2(\ss)$. Thus any $f\in L_2(\ss)$ can be expressed by its
Fourier (or Laplace)
series:$$f=\sum_{\ell=0}^{\infty}{\sum_{k=1}^{Z(d,\ell)} \lb
f,Y_{\ell,k}\rb Y_{\ell,k}},$$
 where
$$\lb f,Y_{\ell,k}\rb =\int_{\ss} f(x)
Y_{\ell,k}(x) \, d\sz(x)$$ are the  Fourier coefficients of $f$.
  The Sobolev space $H^{\sz}
(\ss)$, where $\sz>0$, is defined as the space of all functions in
$L_2(\ss)$ for which the norm
$$\|f\|_{H^\sz}:=\bigg(\sum_{\ell=0}^\infty \Big(1+\ell(\ell+d-1)\Big)^\sz\sum_{k=1}^{N(d,\ell)}\big|\langle
f,Y_{\ell,k}\rangle\big|^2\bigg)^{1/2}$$ is finite. The space $
H^\sz(\ss)$ is a Hilbert space with the inner product $$\langle
f,g\rangle_{ _{H^\sz}}:= \sum_{\ell=0}^\infty
\Big(1+\ell(\ell+d-1)\Big)^\sz\sum_{k=1}^{N(d,\ell)}\langle
f,Y_{\ell,k}\rangle \,\langle g,Y_{\ell,k}\rangle,$$ which induces
the norm $\|\cdot\|_{H^\sz}$. It is easily seen that for $f\in
H^\sz(\ss)$,
\begin{equation}\label{2.1}\|f-S_nf\|_2\le
cn^{-\sz}\|f\|_{H^\sz},\end{equation}where $S_nf$ is given by
\eqref{1.1}.

If  $\sz>d/2$, then  the space $H^\sz(\mathbb{S}^d)$ is
continuously embedded in $C(\mathbb{S}^d)$ and  is  a reproducing
kernel Hilbert space with the reproducing kernel $$
K_\sz(x,y)=\sum_{\ell=0}^\infty
\Big(1+\ell(\ell+d-1)\Big)^{-\sz}\sum_{k=1}^{N(d,\ell)}Y_{\ell,k}(x)Y_{\ell,k}(y).$$

\subsection{Two examples of MZ families} \

We shall give two examples of MZ families  with equal weights. Let
$X$ be a finite subset of $\ss$ with cardinality $d_n:=\dim
\Pi_n^d$.
 The points $\xi_1,\dots, \xi_{d_n}$ in $X$  maximizing a
Vandermonde-type determinant
$$\big|\Delta(\xi_1,\dots,\xi_{d_n})\big|=\big|\det\big(Q_i(\xi_j)\big)_{i,j=1}^{d_n}\big|$$ are called  Fekete points of degree $n$ for $\ss$ (these points are
sometimes called extremal fundamental systems of points, as in
\cite{SW}), where $Q_i, \ i=1,\dots, d_n$ are a basis of
$\Pi_n^d$.  It is known that Fekete points are independent of the
choice of the polynomial basis. For any fixed $\vz>0$, if
$\mathcal X=\{\mathcal X_n\}$ and $\mathcal X_n$ are Fekete points
of degree $\lfloor n(1+\vz)\rfloor$ for $\ss$, then $\mathcal X$
is a MZ family with equal weights $\tau_{n,k}=\frac1{d_n}$ (see
\cite{MO}).

Let $X$ be a finite subset of $\ss$. The mesh norm of $X$ is
$$\rho(X) = \sup_{u\in\ss}\inf_{z\in X} d(u, z),$$where  $d(x,y)$
denotes the standard geodesic distance $\arccos x\cdot y$ between
two points $x$ and $y$ on $\ss$. Let $\mathcal X  =\{\mathcal
X_n\}= \{x_{n,k} : n=1,2,\dots, \, k = 1,\dots, l_n\}$ be a
doubly-indexed set of points in $\ss$. We say that $\mathcal X$ is
uniformly separated if there is a positive number $\vz
> 0$ such that $$d(x_{n,k},x_{n,l})\ge \frac{\vz}{n+1}\ \ {\rm if}\ \   k\neq l$$ for all
$n\in \Bbb N$. Let $\mathcal X=\{\mathcal X_n\}$ be a uniformly
separated array in $\ss$ such that for all $n\ge1$,
$$\rho(\mathcal X_n)\le \frac{\eta}n,
$$ where $\eta<\pi/2$. Then $\mathcal X$ is an MZ family with
the equal weights  $\tau_{n,k}=\frac1{d_n}$ (see \cite{MP}).

However, in \cite{MO} and \cite{MP}, the authors do not give
estimates of the global condition numbers of two MZ  families on
the sphere.

\subsection{ The
weighted  least squares polynomials $ L_n  f$}\

 Now let $\mathcal
X$ be a MZ family with associated weight $\tau$. We use the
weighted discretized inner product
\begin{equation}\label{2.20}\langle
f,g\rangle_{(n)}:=\sum_{k=1}^{l_n}f(x_{n,k})g(x_{n,k})\tau_{n,k}\end{equation}
and the discretized norm $$\|f\|_{(n)}^2=\langle
f,f\rangle_{(n)}.$$

We consider the corresponding orthogonal polynomial projection
$L_n$ onto $\Pi_n^d$ with respect to the weighted discretized
inner product $\langle \cdot,\cdot\rangle_{(n)}$, namely the
weighted  least squares polynomial $ L_n  f$ defined by
$$L_n f=\arg\min_{p\in\Pi_n^d} \|f-p\|_{(n)}^2
=\arg \min_{p\in\Pi_n^d}\
\sum_{k=1}^{l_n}|f(x_{n,k})-p(x_{n,k})|^2\tau_{n,k}.$$ We shall
give  a formal   expression for  $L_nf$. Let $\varphi_i,\
i=1,\dots, d_n$ be an orthonormal basis of $\Pi_n^d$ with respect
to the weighted discrete scalar product \eqref{2.20}, i.e.,
$$\langle
\varphi_i,\varphi_j\rangle_{(n)}=\delta_{i,j}=\Big\{\begin{array}{cl}
1,\ \
   & {\rm if}\ \ i=j,\\
0,\ \     & {\rm if}\ \ i\neq j,
\end{array} \ \ \  1\le i,j\le d_n.$$
Such $\varphi_i,\ i=1,\dots,d_n$ can be obtained by applying
Gram-Schmidt orthogonalization process to the basis $\{Y_{\ell,k}\
|\ k=1,\dots, N(d,\ell),\ \ell=0,1,\dots, n\}$ of $\Pi_n^d$. We
set
$$D_n(x,y)=\sum_{k=1}^{d_n}\varphi_k(x)\varphi_k(y).$$
Clearly, the weighted  least squares polynomial $L_{n}f$ is just
the orthogonal projection of $f$ onto $\Pi_{n}^{d}$  with respect
to the weighted discrete scalar product \eqref{2.20}, i.\! e.,
\begin{equation}\label{2.21}L_{n}f({x})=\langle
f,D_n(\cdot,x)\rangle_{(n)}=\sum_{k=1}^{l_n}
f(x_{n,k})D(x_{n,k},x)\tau_{n,k},\end{equation}and $D_n(x,y)$ is
the reproducing kernel for $\Pi_{n}^{d}$ with respect to the
weighted discrete scalar product \eqref{2.20} satisfying the
following properties:

(1) For any $x,y\in \ss$, $D_{n}(x,y)=D_{n}(y,x)$;

(2) For any fixed $y\in \ss$, $D_{n}(\cdot,y)\in\Pi_{n}^{d}$;

(3) For any $ p\in\Pi_{n}^{d}$ and $x\in \ss$,
\begin{equation}\label{2.22}p(x)=L_np(x)= \langle
p,D_{n}(\cdot,x)\rangle_{(n)}=\sum_{k=1}^{l_n}
p(x_{n,k})D(x_{n,k},x)\tau_{n,k} .\end{equation}

According to the definition of the operator norm and \eqref{2.21},
by the standard argument we have
\begin{equation}\label{2.23}\|L_{n}\|=\max_{x\in \ss}\sum_{k=1}^{l_n}
|D(x_{n,k},x)|\tau_{n,k}.\end{equation}

\subsection{The least squares quadrature}\

  Let
$\mathcal X$ be a MZ family with a weight $\tau$ and let
$E_n(x,y)$ be the reproducing kernel of $\Pi_n^d$. Then the
 frame operator $$T_n p(\cdot) = \sum_{k=1}^{l_n}
\langle p, \tilde e_{n,k}\rangle \tilde
e_{n,k}(\cdot)=\sum_{k=1}^{l_n}\tau_{n,k}p(x_{n,k})E_n(\cdot,x_{n,k})$$
is invertible on $\Pi_n^d$, where $\tilde e_{n,k}=\tau^{ 1/2}_{
n,k}E_n(\cdot,x_{n,k})$. We set $$e_{n,k} = T_n^{-1} (\tilde
e_{n,k})\ \  {\rm and}\ \ w_{n,k} = \langle\tau^{ 1/2}_{
n,k}e_{n,k} , 1\rangle = \tau^{ 1/2}_{ n,k}\int_{\ss} e_{n,k} (x)
d\sz(x).$$ Then the least squares quadrature is defined by
$$I_n(f)=\sum_{k=1}^{l_n}w_{n,k}f(x_{n,k}).$$
In this subsection we shall show  that
\begin{equation}\label{2.24}I_n(f)=\int_{\ss}L_nf(x)d\sz(x)=\sum_{k=1}^{l_n}f(x_{n,k})\tau_{n,k}\int_{\ss}D(x_{n,k},x)d\sz(x).\end{equation}
 In order to prove \eqref{2.24} it suffices to show that $$ \tau^{ 1/2}_{
 n,k} e_{n,k} =\tau_{n,k}D(x_{n,k},\cdot).$$
 Indeed, by \eqref{2.21} and \eqref{2.22} we have for $x,y\in \ss$
 \begin{align*}T_n(D_n(\cdot,y))(x)&=\sum_{k=1}^{l_n}\tau_{n,k}D(x_{n,k},y)E_n(x,x_{n,k})\\&=L_n(E(\cdot,x))(y)=E_n(x,y).\end{align*}It
 follows that
$$\tau^{ 1/2}_{
 n,k} e_{n,k}
 =\tau_{n,k}T_n^{-1}(E_n(\cdot,x_{n,k}))=\tau_{n,k}D(x_{n,k},\cdot).$$Hence,
 \eqref{2.24} holds. By \eqref{2.24} and the H\"older inequality we obtain that
\begin{equation}\label{2.25}\big|\int_{\ss}f(x)d\sz(x)-I_n(f)\big|\le\int_{\ss}|f(x)-L_nf(x)|d\sz(x)\le \|f-L_nf\|_2.\end{equation}

\subsection{Hyperinterpolation and the
weighted  least squares polynomials}\

 Hyperinterpolation
 was originally introduced by I. H. Sloan (see \cite{S}).
It uses the Fourier orthogonal projection of a function which can
be expressed in the form of integrals, but approximates the
integrals used in the expansion by means of a positive cubature
formula. Hence, hyperinterpolation is a discretized orthogonal
projection on polynomial subspaces and provides a polynomial
approximation which relies only on a discrete set of data. In
recent years, hyperinterpolation has attracted much interest, and
a great number of interesting results
 have been obtained (see \cite{CDMV,CDV, CD, D2,  GS,  HAC, HS, DVX, R, R2, S, SW0, SW1, W,  Wang, WHLW, WS, WWW}).

Assume that $Q_n(f)=\sum\limits_{k=1}^{l_n}\tau_{n,k}f(x_{n,k}),\
n=1,2,\dots$ is a sequence of positive quadrature formulas on
$\ss$ which are exact for $\Pi_{2n}^d$, i.e., $\tau_{n,k}>0$, and
for all $f\in \Pi_{2n}^d$,
$$\int_{\ss}f(x)d\sz(x)=Q_n(f)=\sum_{k=1}^{l_n}\tau_{n,k}f(x_{n,k}).$$

For any $p\in\Pi_n^d$, using $f=p^2$ in the above equality we
obtain that  the family $\mathcal X$ is a MZ family with the
constants $A=B=1$ and the global condition number $\kappa$
 is equal to $1$.

 The hyperinterpolation operator $H_n$ on $\ss$
is defined by
\begin{equation*}H_{n}f(x)=\sum\limits_{k=1}^{l_n}\tau_{n,k}f(x_{n,k})E_{n}(x,x_{n,k}),\
\ \ f\in C(\ss).  \end{equation*}

We use the discretized inner product
$$\langle
f,g\rangle_{(n)}:=\sum_{k=1}^{l_n}f(x_{n,k})g(x_{n,k})\tau_{n,k}$$
and the discretized norm $\|f\|_{(n)}^2=\langle f,f\rangle_{(n)}$.
Since the quadrature formula $Q_n$ is exact for $\Pi_{2n}^d$, we
get for all $p,q\in\Pi_n^d$, $$\langle
p,q\rangle=\int_{\ss}p(x)q(x)d\sz(x)=Q_n(p\,q)=\langle
p,q\rangle_{(n)}.$$ Hence,
$$\{Y_{\ell,k}\ |\ k=1,\dots, N(d,\ell),\ \ell=0,1,\dots, n\}$$
forms an orthonormal basis of $\Pi_n^d$ with respect to $\langle
\cdot,\cdot\rangle_{(n)}$. It follows that
$$H_nf(x)=\langle
f,E_n(\cdot,x)\rangle_{(n)}=\sum_{\ell=0}^n\sum_{k=1}^{N(d,\ell)}\langle
f,Y_{\ell,k}\rangle_{(n)}Y_{\ell,k}(x),$$ which means that the
hyperinterpolation operator $H_n$ is just the discretized
orthogonal projection on the polynomial subspace $\Pi_n^d$ with
respect to  $\langle \cdot,\cdot\rangle_{(n)}$, and satisfies
$$\|f-H_nf\|_{(n)}^2=\min_{p\in\Pi_n^d}\|f-p\|_{(n)}^2=\min_{p\in\Pi_n^d}\sum_{k=1}^{l_n}|f(x_{n,k})-p(x_{n,k})|^2\tau_{n,k}.$$

Hence, the hyperinterpolation $H_nf$ is just  the weighted least
squares polynomial for a sequence of positive quadrature formulas
on $\ss$ which are exact for $\Pi_{2n}^d$.

On the other hand, if the global condition number $\kappa$ of a MZ
family $\mathcal X$ with the weight $\tau$ is equal to $1$, then
for all $p\in \Pi_n^d$,  we have
$$A\|p\|_2^2=\sum_{k=1}^{l_n}|p(x_{n,k})|^2\tau_{n,k}.$$It follows
that for $f,g\in\Pi_n^d$,
$$A\langle
f,g\rangle=\sum_{k=1}^{l_n}f(x_{n,k})g(x_{n,k})\tau_{n,k}.$$ Since
any function in $\Pi_{2n}^d$ can be expressed as a linear
combination of the product of the functions $x^\az$ and $x^\beta$,
we obtain that $\mathcal X$ determines a sequence of positive
quadrature formulas $Q_n(f)=\frac
1A\sum\limits_{k=1}^{l_n}\tau_{n,k}f(x_{n,k})$ on $\ss$ which are
exact for $\Pi_{2n}^d$,  where, $x\in\ss,\ \az,\beta\in\Bbb
N_0^{d+1}, |\az|\le n,\ |\beta|\le n,\ x^\az:=x_1^{\az_1}\dots
x_{d+1}^{\az_{d+1}}, \ |\az|=\az_1+\dots+\az_{d+1}$. This means
that the weighted least squares operators $L_n$ are just the
hyperinterpolation operators $H_n$.

Hence,  the hyperinterpolation
 is just  the weighted least squares polynomial for a   MZ
family with the global condition number $\kappa=1$, and the
weighted least squares polynomial may be viewed as a
generalization of the hyperinterpolation.

\section{Proofs of Theorems 1.2 and 1.4}

In order to prove Theorem 1.2, we shall use the following lemma.

\begin{lem}\label{3-2-lem}\cite[Theorem 2.1]{D2}
Let  $\Ga$ be   a finite subset of  $\sph$, and  let  $\{
\mu_\omega:\ \ \omega\in \Ga\}$ be   a set of positive  numbers
satisfying
\begin{equation*} \sum_{ \omega\in\Ga} \mu_\omega
|f(\omega)|^{p_0} \leq C_1 \int_{\sph} |f(x)|^{p_0}\, d\sa(x),\ \
\forall f\in \Pi_{N}^d,
\end{equation*} for some $0<p_0<\infty$ and some positive integer $N$. If
 $0< q<\infty$, $M\ge
N$ and $f\in \Pi_M^d$, then
\begin{equation*} \sum_{ \omega\in \Ga} \mu_\omega
|f(\omega)|^q \leq C C_1 \( \f MN\)^{d} \int_{\sph} |f(y)|^q \,
d\sa(y),\end{equation*} where  $C>0$ depends only on $d$ and  $q$.
\end{lem}

\noindent{\it Proof of Theorem 1.2}.

 We rewrite the orthonormal
basis $\{Y_{\ell,k}\ |\ k=1,\dots, N(d,\ell),\ \ell=0,1,\dots,
n\}$ of $\Pi_n^d$ as $\{\phi_k\ |\ k=1,\dots,d_n\}$. Let $${\bf
y}_n=(\tau_{n,1}^{1/2}f(x_{n,1}),\dots,
\tau_{n,l_n}^{1/2}f(x_{n,l_n}))^*$$ be the given sampling data
column vector, $U_n$ be the $l_n\times d_n$ matrix with entries
$$(U_n)_{kl}=\tau_{n,k}^{1/2}\phi_l(x_{n,k}), \ \ \
k=1,\dots,l_n,\ l=1,\dots,d_n,$$ and let $R_n=U_n^*U_n$, where
$U_n^*$ is the conjugate transpose of the matrix $U_n$. Suppose
that the solution to the weighted least squares problem
\eqref{1.3} is the polynomial
$$L_nf=\sum_{k=1}^{d_n} a_{n,k}\phi_k\in\Pi_n^d$$with coefficient
column vector ${\bf a}_n$. According to  the standard formula for
the solution of a least squares problem by means of the
Moore-Penrose inverse $U_n^+=(U_n^*U_n)^{-1}U_n^*=R_n^{-1}U_n^*$,
we obtain that
 $${\bf a}_n=U_n^+{\bf y}_n= R_n^{-1}U^*_n{\bf y}_n.$$
For any $c\in \Bbb R^{d_n}$, $p=\sum_{k=1}^{d_n}c_k\phi_k$, we
have
$$\|p\|_2^2=|c|^2\ \ {\rm and}\ \ (R_nc,c)=(U_nc,U_nc)=\sum_{k=1}^n |p(x_{n,k}|^2\tau_{n,k}. $$
It follows from \eqref{1.2} that
$$ A|c|^2\le (R_nc,c)\le B|c|^2, $$
which means that the spectrum of every $R_n$ is contained in the
interval $[A, B]$. We obtain  that the operator norm of $R_n^{-1}$
is bounded by $A^{-1}$ and the operator norm of $U_n^*$ is bounded
by
$$\|U_n^*\|=\|U_n\|=\|U_n^*U_n\|^{1/2}=\|R_n\|^{1/2}\le
B^{1/2},$$where the operator norm  of a matrix $A$ is defined by
$\|A\|=\sup_{|x|=1}|Ax|.$

We use the orthogonal decomposition
$$\|f-L_nf\|_2^2=\|f-S_nf\|_2^2+\|S_nf-L_nf\|_2^2.$$ We recall that
$$S_nf=\sum_{k=1}^{d_n}f_{n,k}\phi_k$$with coefficient column
vector ${\bf f}_n$, where $f_{n,k}=\langle f,\phi_k\rangle$. By
the Parseval equality we obtain that
\begin{align*}\|S_nf-L_nf\|_2^2&=|{\bf f}_n-{\bf a}_n|^2=|{\bf
f}_n-R_n^{-1}U_n^*{\bf y}_n|^2\\ &=|R_n^{-1}U_n^*(U_n{\bf
f}_n-{\bf y}_n)|^2\\ &\le A^{-2}B|U_n{\bf f}_n-{\bf
y}_n|^2.\end{align*} Finally, we have
$$(U_n{\bf f}_n)_k=\tau_{n,k}^{1/2}S_nf(x_{n,k}).$$Hence,
$$|U_n{\bf f}_n-{\bf
y}_n|^2=\sum_{k=1}^{l_n}|f(x_{n,k})-S_nf(x_{n,k})|^2\tau_{n,k}=\|f-S_nf\|_{(n)}^2.$$
It follows that \begin{equation}\label{3.1} \|S_nf-L_nf\|_2^2\le
A^{-2}B \|f-S_nf\|_{(n)}^2
\end{equation}

For $f\in H^\sz(\ss),\ \sz>d/2$, we define
$$A_0f=S_{n}f, \ \ A_kf=S_{2^{k}n}f-S_{2^{k-1}n}f\ \
{\rm for}\ \ k\geq 1.$$ Then for $f\in H^\sz(\ss),\  s>d/2$, the
series $\sum\limits_{k=0}^\infty A_k f(x)$ converges   to $f(x)$
uniformly in $\ss$. By \eqref{2.1} we have
\begin{equation}\label{3.2} \|A_kf\|_2\le \|f-S_{2^kn}f\|_2+  \|f-S_{2^{k-1}n}f\|_2\ll
2^{-k\sz}n^{-\sz}\|f\|_{H^\sz}.\end{equation} Note that
$$f(x)-S_nf(x)=\sum_{k=1}^\infty A_kf(x),$$and the right series
converges uniformly on $\ss$. Hence, by \eqref{3.1} and the
trigonometric  inequality  we have
$$\|S_nf-L_nf\|_2\le A^{-1}B^{1/2}\|f-S_nf\|_{(n)}\le A^{-1}B^{1/2}\sum_{k=1}^\infty \|A_kf\|_{(n)}.$$
Note that $A_kf\in\Pi_{2^kn}^d$. Using Lemma 3.1 with $p_0=q=2$
and \eqref{3.2}, we obtain that
\begin{align*}\|A_kf\|_{(n)}^2&=\sum_{k=1}^{l_n}|A_kf(x_{n,k})|^2\tau_{n,k}\\ &\le
cB\Big(\frac{2^kn}{n}\Big)^d\int_{\ss}|A_kf(x)|^2\,d\sz(x)\\
&\ll B2^{kd}2^{-2k\sz}n^{-2\sz}\|f\|_{H^\sz}^2.\end{align*} We
have
\begin{align*}\|S_nf-L_nf\|_2&\le A^{-1}B^{1/2} \sum_{k=1}^\infty
\|A_kf\|_{(n)}\\
&\ll A^{-1}B\sum_{k=1}^\infty 2^{kd/2}2^{-k\sz}n^{-\sz}\|f\|_{H^\sz}\\
&\ll A^{-1}Bn^{-\sz}\|f\|_{H^\sz}.\end{align*} It follows that
\begin{align*}\|f-L_nf\|_2&=\big(\|f-S_nf\|_2^2+\|S_nf-L_nf\|_2^2\big)^{1/2}\\ &\ll
(1+\kappa^2)^{1/2}n^{-\sz}\|f\|_{H^\sz}.
\end{align*}
This completes the proof of \eqref{1.4}. Inequality \eqref{1.5}
follows from \eqref{1.4} and \eqref{2.25} directly.

 The proof of Theorem 1.2 is now finished. $\hfill\Box$

\begin{rem}The inequality \eqref{3.1} was proved in \cite{Gr}. For the convenience of the reader, we provide details of the proof.\end{rem}

\noindent{\it Proof of Theorem 1.4}.

In order to give the lower estimate of the operator norm
$\|L_n\|$, we use the  Daugavet theorem to obtain that
$$\|L_n\|\ge \|S_n\|\asymp n^{(d-1)/2}.$$
 So it
suffices to get the upper estimate of $\|L_n\|$.

Let the functions
 $$\Phi_n(x)=\inf_{p\in\Pi_n^d,\ |p(x)|=1} \|f\|_2^2,\ \ x\in\ss,$$ and $$\Psi_n(x)=\inf_{p\in\Pi_n^d,\ |p(x)|=1}
 \|f\|_{(n)}^2,\ \ x\in\ss$$ be the Christoffel functions with respect to the
 inner $\langle \cdot,\cdot\rangle$ and $\langle
 \cdot,\cdot\rangle_{(n)}$, respectively. We remark that estimates for Christoffel functions are useful in comparing different norms of functions in $\Pi_n^d$, and they are also basic tools in the theory of orthogonal polynomials.

It
 follows from \eqref{1.2} that
 \begin{equation}\label{3.10}A\Phi_n(x)\le \Psi_n(x)\le
 B\Phi_n(x),\ \ x\in\ss.\end{equation}
 Now  let $E_n(x,y)$ and $D_n(x,y)$
be the reproducing kernels of $\Pi_n^d$ with respect to the
 inner $\langle \cdot,\cdot\rangle$ and $\langle
 \cdot,\cdot\rangle_{(n)}$, respectively.  According to \cite[Theorem 3.5.6]{DX}, we have
 $$E_n(x,x)^{-1}=\Phi_n(x) \ \ {\rm and}\  \ D_n(x,x)^{-1}=\Psi_n(x).$$It
 follows from \eqref{3.10} that \begin{equation}\label{3.11}B^{-1}E_n(x,x)\le D_n(x,x)\le A^{-1}
 E_n(x,x), \ \ x\in\ss.\end{equation} Using $p=1$ in \eqref{1.2} we get that
\begin{equation}\label{3.12}A\le \sum_{k=1}^{l_n}\tau_{n,k}\le
B.\end{equation}
 Using \eqref{2.23}, the Cauchy-Schwartz inequality, \eqref{2.22}, \eqref{3.12}, \eqref{3.11}, and \eqref{2.01},    we have
 \begin{align*}\|L_n\|&=\max_{x\in\ss}\sum_{k=1}^{l_n}\tau_{n,k}|D_n(x,x_{n,k})|\\&\le
 \max_{x\in\ss}\Big(
 \sum_{k=1}^{l_n}\tau_{n,k}(D_n(x,x_{n,k}))^2\Big)^{1/2}
 \Big(\sum_{k=1}^{l_n}\tau_{n,k} \Big)^{1/2}\\ &\le B^{1/2}\max_{x\in\ss}
 D_n(x,x)^{1/2}\\ &\ll\kappa^{1/2} \max_{x\in\ss} E_n(x,x)^{1/2}\\ &\ll
 \kappa^{1/2}
 n^{d/2}.\end{align*}

This completes the proof of  Theorem 1.4. $\hfill\Box$

\end{document}